\documentclass[12pt]{amsart}
\usepackage[cp1251]{inputenc}
\usepackage[russian]{babel}
\usepackage{amssymb}

\theoremstyle{plain}

\newtheorem{Statement}{Утверждение}
\newtheorem{Theorem}{Теорема}

\newcommand{\intl}{\mathop{\int}\limits}

\newcommand{\supl}{\mathop{\sup}\limits}
\renewcommand{\Re}{\operatorname{Re}}
\renewcommand{\Im}{\operatorname{Im}}

\renewcommand{\le}{\leqslant}
\renewcommand{\ge}{\geqslant}

\def\al{\alpha}
\def\la{\lambda}

\def\w{\operatorname{w}}
\def\v{\operatorname{v}}
\def\ur{u_{\mathcal{R}}}
\def\ui{u_{\mathcal{I}}}

\def\N{{\mathbb N}}

\def\C{{\mathbb C}}

\newlength{\lenun}
\newlength{\lendu}
\begin{document}
{\Large {\bf \centerline{О собственных функциях оператора
Штурма--Лиувилля}} {\bf \centerline{с потенциалами из пространств
Соболева}}}\footnote{Работа поддержана грантом РФФИ No.
09-01-90408}

\bigskip

{\centerline{Савчук А.~М.}}

\bigskip
В работе изучается оператор Штурма--Лиувилля $L$ в пространстве $L_2[0,\pi]$, порождённый дифференциальным выражением $l(y)=-y''+q(x)y$ и краевыми условиями Дирихле $y(0)=y(\pi)=0$. Потенциал $q\in W_2^{-1}[0,\pi]$ предполагается комплекснозначным (здесь $W_2^{-1}[0,\pi]$ --- пространство Соболева). Операторы такого вида были определены в работе \cite{SS1} и подробно изучены в работах \cite{SS1}--\cite{SS3}. В частности, было доказано, что спектр $L$ чисто дискретен и были получены асимптотические формулы для собственных значений и собственных функций. При этом в формулах для собственных значений были выписаны два асимптотических члена, а для собственных функций был найден только главный член асимптотики. Цель настоящей работы --- получить второй член в асимптотике собственных функций. Интерес к такого рода оценкам вызван недавними работами И.В.Садовничей \cite{Sa1}, где эти результаты используются в теоремах о равносходимости.

Нам потребуются некоторые обозначения и утверждения, доказанные в работе \cite{SS2}. Мы обозначим через $u(x)$ --- обобщённую первообразную функции $q$ (она определена с точностью до константы), $\ur=\Re u$, $\ui=\Im u$. Обозначим также
\begin{gather*}
b(x,\la):=\intl_0^x u(t)\sin(2\la^{1/2}t)dt, \quad
a(x,\la):=\intl_0^x u(t)\cos(2\la^{1/2}t)dt,\\
B(x,\la):=\intl_0^x u^2(t)\sin(2\la^{1/2}t)dt,\quad
A(x,\la):=\intl_0^x u^2(t)\cos(2\la^{1/2}t)dt,\\
U(x):=\intl_0^x u^2(t)dt,\quad
\w(x,\la):=\intl_0^x\intl_0^tu(t)u(s)\cos(2\la^{1/2}t)\sin(2\la^{1/2}s)dsdt,\\
\upsilon(x,\la):=b(x,\la)+\frac12\la^{-1/2}U(x)+2\w(x,\la)-\frac12\la^{-1/2}A(x,\la),\\
\Upsilon(\la):=\sup\limits_{0\le x\le\pi}\left(|b(x,\la)|+|a(x,\la)|+2|\w(x,\la)|+\frac12|\la^{-1/2}A(x,\la)|\right)+|\la|^{-1/2}\|u\|^2_{L_2},\\
a(x,n^2)=a_{2n}(x),\quad b(x,n^2)=b_{2n}(x),\\
A(x,n^2)=A_{2n}(x),\quad B(x,n^2)=B_{2n}(x),\\
\w(x,n^2)=\w_{2n}(x),\quad \upsilon(x,n^2)=\upsilon_{2n}(x),\\
b_n:=\intl_0^\pi u(t)\sin(nt)dt,\quad a_n:=\intl_0^\pi u(t)\cos(nt)dt,
\end{gather*}
\begin{gather*}
A_n:=\intl_0^\pi u^2(t)\cos(nt)dt,\quad B_n:=\intl_0^\pi u^2(t)\sin(nt)dt,\\
\w_n:=\intl_0^\pi\intl_0^t u(t)u(s)\cos(2nt)\sin(2ns)dsdt,\quad
U(\pi)=\intl_0^\pi u^2(t)dt,\\
\mu_n:=-\frac1\pi \upsilon(\pi,n^2)=-\frac1\pi\upsilon_{2n}(\pi)= -\frac1{\pi}\intl_0^\pi u(t)\sin(2nt)dt +\frac1{2\pi n}\intl_0^\pi u^2(t)\cos(2nt)dt-\\
-\frac2{\pi}\intl_0^\pi\intl_0^t u(t)u(s)\cos(2nt)\sin(2ns)dsdt -\frac1{2\pi n}\intl_0^\pi  u^2(t)dt=\\
=-\frac1\pi b_{2n}+\frac1{2\pi n}A_{2n}-\frac2\pi\w_{2n}-\frac1{2\pi n}U(\pi).
\end{gather*}

Далее мы будем проводить оценки в комплексной $\la$--плоскости внутри областей, ограниченных параболами
$$
P_\al=\{\la\in\C\,\vert\ \Re\la>1,\ |\Im\sqrt\la|<\al\},
$$
а в $z$--плоскости ($z=\sqrt\la$) внутри полуполос
$$
\Pi_\al=\{z\in\C\,\vert\ \Re z>1,\ |\Im z|<\al\}.
$$
Всюду далее, рассматривая функцию $z=\sqrt{\la}$, подразумеваем выбор
её главной ветви, принимающей положительные значения при $\la>0$.

Введённые выше функции будут участвовать в асимптотиках решений, причём $\Upsilon^2(\la)$ будет служить для оценки остатков. Функцию $\Upsilon(\la)$ можно заменить более простым выражением
$$
\Upsilon_1(\la)=\sup\limits_{0\le x\le\pi}\left(|b(x,\la)|+|a(x,\la)|\right)+|\la|^{-1/2}\|u\|^2_{L_2},
$$
так как внутри парабол $P_\al$ легко получить оценку
$$
\Upsilon_1(\la)\le\Upsilon(\la)\le M\Upsilon_1(\la),
$$
где $M=2\sqrt\pi(1+\|u\|_{L_2})\rm ch(2\pi\al)$.

Представим основные результаты о спектре оператора $L$.
\begin{Statement}(см. \cite[Теорема 1.5]{SS2})\label{st:1}
Оператор \(L\) имеет непустое резольвентное множество и спектр его дискретен.
\end{Statement}
Обозначим через $\omega(x,\la)$ решение
дифференциального уравнения $-\omega''+q\omega=\la\omega$ с
начальными условиями $\omega(0,\la)=0$, $\omega^{[1]}(0,\la)=1$
(здесь $\omega^{[1]}=\omega'-u\omega$ --- первая
квазипроизводная). Ясно, что нули целой функции $\omega(\pi,\la)$
совпадают с собственными значениями оператора $L$. В силу теоремы
существования и единственности, геометрическая кратность каждого
собственного значения равна $1$. В случае вещественного потенциала
присоединенные функции отсутствуют, т.е. все собственные значения
являются простыми, но в общем случае это не так. Обозначим через
$\{\mu_k\}_1^\infty$ все нули функции $\omega(\pi,\la)$ \emph{без
учета кратности} (т.е. $\mu_j=\mu_k$ только при $j=k$), причем
нумерацию будем вести в порядке возрастания модуля, а в случае
совпадения модулей --- по возрастанию аргумента, значения которого
выбираются из полуинтервала $(-\pi,\pi]$. Цепочкой из собственной
и присоединенных функций, отвечающей собственному значению
$\mu_k$, называют систему функций
$\{y_k^0,\,y_k^1,\,\dots,\,y_k^{p_k-1}\}$, где
$(L-\mu_kI)y_k^0=0$, $(L-\mu_kI)y_k^j=y_k^{j-1}$ для всех
$j=1,\,\dots,\,p_k-1$. Присоединенные функции определены
неоднозначно и, следуя работе \cite{Ke}, \emph{канонической
цепочкой} назовем любую из цепочек, имеющих максимальную длину.
Длину канонической цепочки (число $p_k$) назовем
\emph{алгебраической кратностью} собственного значения $\mu_k$. В
нашем случае каноническую цепочку можно предъявить явно, что мы
сейчас и проделаем.

Пусть $\mu_k$ --- ноль функции $\omega(\pi,\la)$ кратности $q_k$.
Заметим, что функции $\omega^{(j)}_{\la}(x,\la)$,
$j=1,2,\dots,q_k-1$ удовлетворяют дифференциальным уравнениям
$l(\omega^{(j)}_{\la})=\la\omega^{(j)}_{\la}+\omega^{(j-1)}_\la$,
причем $\omega^{(j)}_{\la}(0,\la)\equiv0$. Кроме того,
$\omega^{(j)}_\la(\pi,\mu_k)=0$, поскольку $\mu_k$ есть ноль
кратности $q_k$. Таким образом, функции
$\omega^{(j)}_\la(x,\mu_k)$, $j=0,1,\dots,q_k-1$ образуют цепочку
из собственной и присоединенных функций, отвечающую собственному
значению $\mu_k$. Как и любая цепочка из собственной и
присоединенных функций, эта система линейно независима, а значит,
порождает подпространство размерности $q_k$. Теперь остается
отметить (см. \cite[Гл. I, п.3]{Na}), что эта цепочка является
канонической, т.е. имеет максимальную длину ($p_k=q_k$).

Нам будет удобно ввести и другую нумерацию собственных и
присоединенных функций.  Под системой собственных и присоединенных
функций $\{y_n(x)\}_{n=1}^\infty$ мы будем везде далее понимать
систему
$\bigcup\limits_{k=1}^{\infty}\left\{\frac{\omega^{(j)}_\la(\cdot,\,\mu_k)}{\|\omega(\cdot,\,\mu_k)
\|_{L_2}}\right\}_{j=0}^{q_k-1}$ (таким образом,
$\|y_n\|_{L_2}=1$, если $y_n$ --- собственная функция). Через
$\{\la_n\}_{n=1}^\infty$ обозначим нули функции $\omega(\pi,\la)$
в порядке возрастания модуля \emph{с учетом кратности} (тогда
$Ly_n=\la_ny_n$ для любого $n\in\N$).

Построим теперь биортогональную систему $\{\v_n\}_{n=1}^\infty$ к
системе $\{y_n\}_{n=1}^\infty$. Прежде всего заметим, что при
любом $k$ система $\{z_k^j\}_{j=0}^{q_k-1}$, где
$z_k^j(x)=\overline{\omega^{j}_\la(x,\mu_k)}$, является системой
из собственной и присоединенных функций для оператора
$L^*=-\frac{d^2}{dx^2}+\overline{q}$, отвечающей собственному
значению $\overline{\mu_k}$. Покажем теперь, как с помощью
конечных линейных комбинаций системы
$\bigcup\limits_{k=1}^{\infty}\{z_k^j\}_{j=0}^{q_k-1}$ построить
биортогональную к  $\{y_n\}_1^\infty$ систему. Заметим, что если
функция $z$ лежит в корневом подпространстве оператора $L^*$,
отвечающем собственному значению $\overline{\mu_l}$, а функция $y$
--- в корневом подпространстве оператора $L$, отвечающем
собственному значению $\mu_k$, где $k\ne l$, то функции $y$ и $z$
ортогональны. Действительно, для собственных функций
$\mu_k(y,z)=(Ly,z)=(y,L^*z)=\mu_l(y,z)$, т.е. $(y,z)=0$. Если
теперь $y$ --- первая присоединенная функция, а $z$ ---
собственная функция, то, учитывая, что ортогональность собственных
функций уже доказана, $\mu_k(y,z)=(Ly,z)=(y,L^*z)=\mu_l(y,z)$.
Дальнейшее очевидно. Таким образом, достаточно построить
биортогональную систему в каждом корневом подпространстве по
отдельности. В случае простого собственного значения это легко:
$\v_n(x)=\overline{y_n(x)}/(y_n(x),\overline{y_n(x)})$. Знаменатель
здесь отличен от нуля, так как интегрированием по частям легко
получить равенство
$\intl_0^{\pi}\omega^2(x,\mu_k)dx=\omega'_{\lambda}(\pi, \mu_k)
\omega'_x(\pi, \mu_k)$. Для собственных значений кратности $p_k>1$
мы сошлемся на \cite[Гл. 1, \S 3]{Mar}:
\begin{Statement}
$$
(z_k^j,\omega_k^i)=\begin{cases} 0, \quad\quad &i+j<p_k-1, \\
a^{i+j}_k, \quad &i+j\ge p_k-1,
\end{cases}
$$
причем $a^{p_k-1}_k\neq 0$.
\end{Statement}

Тогда элементы биортогональной системы имеют вид
$$
\v^{p_k-j}_k(x)=b^{j-1}_kz^0_k(x)+b^{j-2}_kz^1_k(x)+\dots
+b^0_kz^{j-1}_k(x),
$$
где $b^0_k=\dfrac{1}{a^{p_k-1}_k}$,
$b^1_k=-\dfrac{b^0_ka^{p_k}_k}{a^{p_k-1}_k}$,
$b^2_k=-\dfrac{b^1_ka^{p_k}_k+b^0_ka^{p_k+1}_k}{a^{p_k-1}_k}$,
$\dots$,
$$b^{p_k-1}_k=-\dfrac{b^{p_k-2}_ka^{p_k}_k+b^{p_k-3}_ka^{p_k+1}_k+\dots
+b^0_ka^{2p_k-2}_k}{a^{p_k-1}_k}.$$

Следуя работе \cite{SS2} функцию $\omega(x,\la)$ будем искать в виде $\omega(x,\la)=r(x,\la)\sin\theta(x,\la)$ .

\begin{Statement}Функции $\theta(x,\la)$ и $r(x,\la)$ удовлетворяют системе уравнений
\begin{equation}\label{thr}
\left\{
\begin{array}{l}
\theta'(x,\la)=\la^{\frac1{2}}+\la^{-\frac1{2}}u^2(x)\sin^2\theta(x,\la)+u(x)\sin2\theta(x,\la),\\
r'(x,\la)=-r(x,\la)\left[u(x)\cos 2\theta(x,\la)+\frac12\la^{-1/2}u^2(x)\sin 2\theta(x,\la)\right].
\end{array}\right.
\end{equation}
\end{Statement}

\begin{Statement}(см. \cite[Лемма 2.1]{SS2})\label{st:3} Пусть $\al>0$ --- произвольное фиксированное число, $P_\al$ --- область, ограниченная параболой $|\Im\sqrt{\la}|<\al$. Тогда существует число $\mu$ (зависящее только от функции $u$ и $\al$) такое, что при любых $\la\in P_\al$, $\Re\la>\mu$, первое уравнение системы \eqref{thr} имеет единственное решение $\theta(x,\la)$, определённое при всех $0\le x\le\pi$ и удовлетворяющее начальному условию $\theta(0,\la)=0$. Это решение допускает представление
\begin{equation}
\theta(x,\la)=\la^{1/2}x+b(x,\la)+\frac12\la^{-1/2}U(x)+2\w(x,\la)-\frac12\la^{- 1/2}A(x,\la)+\rho(x,\la),
\label{asth}
\end{equation}
где
$$
|\rho(x,\la)|\le M\Upsilon^2(\la),\quad \la\in P_\la,\quad \Re\la>\mu,
$$
причём выбор постоянной $M$ зависит от $u$ и $\al$, но не зависит от $x$, $\la$.
\end{Statement}

\begin{Statement}(см. \cite[Лемма 2.2]{SS2})\label{st:4} Пусть $\al>0$ --- произвольное фиксированное число, а $P_\al$ --- область, ограниченная параболой $|\Im\sqrt{\la}|<\al$. Пусть $\theta(x,\la)$ --- решение первого уравнения системы \eqref{thr} с начальным условием $\theta(0,\la)=0$.  Тогда решение $r(x,\la)$ второго уравнения системы \eqref{thr} с
начальным условием $r(0,\la)=1$ допускает представление
\begin{equation}
r(x,\la)=1-a(x,\la)-\frac12\la^{-1/2}B(x,\la)+\rho(x,\la)\quad\la\in P_\al,\ \Re\la>\mu,
\label{asr}
\end{equation}
где
$$
|\rho(x,\la)|\le M\Upsilon^2(\la).
$$
Здесь $\mu$ и $M$ --- числа, зависящее только от $\al$ и $u(x)$.
\end{Statement}

Представим основные результаты о спектре оператора $L$.

\begin{Statement}(см. \cite[Теорема 3.1]{SS2})\label{st:5} Пусть $\{\la_n\}_1^\infty$ --- собственные значения оператора $L$, занумерованные в порядке возрастания модуля с учётом алгебраической кратности. Тогда
\begin{equation}\label{lan}
\la^{1/2}_n=n+\mu_n+\rho_n, \quad\text{где}\ |\rho_n|\le M\Upsilon^2(\la_n),
\end{equation}
причем выбор постоянной $M$ зависит от $u$, но не зависит от $n$.
\end{Statement}

\begin{Statement}(см. \cite[Лемма 2.3]{SS2}\label{st:6} Пусть несгущающаяся последовательность комплексных чисел
$\{z_n\}$ лежит в полосе $|\Im z|<\al$, где $\al>0$ --- произвольное
число. Пусть оператор $T_x:\,L_2[0,\pi]\to l_2$ определён равенством
$$
T_xf=\{c_n\}_{n=1}^\infty,\quad\text{где}\ c_n=\intl_0^xf(t)e^{iz_nt}dt.
$$
Тогда $\|T_x\|\le M$ для некоторой постоянной $M$, зависящей от последовательности $\{z_n\}$, но не зависящей от $x\in[0,\pi]$.
\end{Statement}
Таким образом, согласно утверждениям 5 -- 7, собственные значения $\{\la_n\}_1^\infty$ оператора $L$ подчинены асимптотике
\begin{equation}\label{evres}
\la_n^{1/2}=n+\mu_n+\rho_n,\quad\text{где }\{\rho_n\}_1^\infty\in l_1.
\end{equation}

Перейдём к результатам о системе собственных и присоединённых функций.
\begin{Theorem}
Начиная с некоторого номера все собственные значения однократны, а нормированные собственные функции имеют вид
\begin{gather*}
y_n(x)=\frac{\omega(x,\la_n)}{\|\omega(x,\la_n)\|}=
\sqrt{\frac2\pi}\sin(nx)\left(1+\frac1\pi\intl_0^\pi(\pi-t)\ur(t)\cos(2nt)dt+\right.\\
\left.+\frac1{2\pi n}\intl_0^\pi(\pi-t)(\ur^2(t)-\ui^2(t))\sin(2nt)dt-
\intl_0^xu(t)\cos(2nt)dt-\frac1{2n}\intl_0^xu^2(t)\sin(2nt)dt\right)+
\end{gather*}
\begin{gather}\label{asef}
+\sqrt{\frac2\pi}x\cos(nx)\left(-\frac1\pi\intl_0^\pi u(t)\sin(2nt)dt+\frac1{2\pi n}\intl_0^\pi u^2(t)\cos(2nt)dt-\right.\\
\left.-\frac2{\pi}\intl_0^\pi\intl_0^tu(t)u(s)\cos(2nt)\sin(2ns)dsdt-\frac1{2\pi n}\intl_0^\pi u^2(t)dt\right)+\notag
\end{gather}
\begin{gather*}
+\sqrt{\frac2\pi}\cos(nx)\left(\intl_0^xu(t)\sin(2nt)dt+\frac1{2n}\intl_0^xu^2(t)dt-
\frac1{2n}\intl_0^xu^2(t)\cos(2nt)dt+\right.\\
\left.+2\intl_0^x\intl_0^tu(t)u(s)\cos(2nt)\sin(2ns)dsdt\right)+\rho_n(x),\quad\text{где }
\{\|\rho_n(x)\|_{C[0,\pi]}\}_1^\infty\in l_1.
\end{gather*}
Соответствующие функции биортогональной системы имеют вид
\begin{gather}\label{asbiort}
\v_n(x)
=\sqrt{\frac2\pi}\sin(nx)\left(1+\frac1\pi\intl_0^\pi(\pi-t)(\ur(t)+2i\ui(t))\cos(2nt)dt-
\intl_0^x\overline{u}(t)\cos(2nt)dt+\right.\notag\\
\left.+\frac1{2\pi n}\intl_0^\pi(\pi-t)(\ur^2(t)-\ui^2(t)+4i\ur(t)\ui(t))\sin(2nt)dt
-\frac1{2n}\intl_0^x\overline{u}^2(t)\sin(2nt)dt\right)+\notag\\
+\sqrt{\frac2\pi}x\cos(nx)\left(-\frac1\pi\intl_0^\pi\overline{u}(t)\sin(2nt)dt+\frac1{2\pi n}\intl_0^\pi \overline{u}^2(t)\cos(2nt)dt-\right.\notag\\
\left.-\frac2{\pi}\intl_0^\pi\intl_0^t\overline{u}(t)\overline{u}(s)\cos(2nt)\sin(2ns)dsdt-\frac1{2\pi n}\intl_0^\pi \overline{u}^2(t)dt\right)+\\
+\sqrt{\frac2\pi}\cos(nx)\left(\intl_0^x\overline{u}(t)\sin(2nt)dt+
\frac1{2n}\intl_0^x\overline{u}^2(t)dt+2\intl_0^x\intl_0^t\overline{u}(t)\overline{u}(s)\cos(2nt)\sin(2ns)dsdt-
\right.\notag\\\left.-
\frac1{2n}\intl_0^x\overline{u}^2(t)\cos(2nt)dt\right)+\widetilde{\rho}_n(x),\quad
\text{где }\{\|\widetilde{\rho}_n(x)\|_{C[0,\pi]}\}_1^\infty\in l_1.\notag
\end{gather}

\end{Theorem}
\begin{proof}
Для удобства будем обозначать через $\rho(x,\la)$ любую функцию, удовлетворяющую оценке $\supl_{x\in[0,\pi]}|\rho(x,\la)|<M\Upsilon^2(\la)$, где $M$ не зависит от $\la$, а $\la\to\infty$ внутри некоторой полуполосы $\Pi_\al$. Будем обозначать через $\rho_n$ любую последовательность из пространства $l_1$, а через $\rho_n(x)$ --- любую функцию, для которой $\left\{\|\rho_n(x)\|_{C[0,\pi]}\right\}_1^\infty\in l_1$. Тогда в силу утверждений 2 -- 4, для любого $\la\in\Pi_\al$ имеем
\begin{equation}\label{omega}
\omega(x,\la)=\sin(\la^{1/2}x)\left(1-a(x,\la)-\frac12\la^{-1/2}B(x,\la)\right)+
\cos(\la^{1/2}x)\upsilon(x,\la)+\rho(x,\la).
\end{equation}
Далее, согласно \eqref{evres}, $\la_n^{-1/2}=n^{-1}+\rho_n$,
$$
\sin(\la^{1/2}_n\xi)=\sin(n\xi)+\mu_n\xi\cos(n\xi)+\rho_n(\xi),\quad
\cos(\la^{1/2}_n\xi)=\cos(n\xi)-\mu_n\xi\sin(n\xi)+\rho_n(\xi).
$$
Тогда
\begin{gather*}
b(x,\la_n)=\intl_0^x u(t)\sin(2\la_n^{1/2}t)dt=b_{2n}(x)+\mu_n a_{2n}(x)=b_{2n}(x)+\rho_n(x);\\
a(x,\la_n)=\intl_0^x u(t)\cos(2\la^{1/2}_nt)dt=a_{2n}(x)-\mu_n b_{2n}(x)=a_{2n}(x)+\rho_n(x);
\end{gather*}
\begin{gather*}
\frac12\la_n^{-1/2}\intl_0^xu^2(t)\sin(2\la_n^{1/2}t)dt=\frac1{2n}B_{2n}(x)+\rho_n(x);\\
\frac12\la_n^{-1/2}\intl_0^xu^2(t)\cos(2\la_n^{1/2}t)dt=\frac1{2n}A_{2n}(x)+\rho_n(x);\\
\w(x,\la_n)=\intl_0^x\intl_0^tu(t)u(s)\sin(2\la_n^{1/2}s)\cos(2\la_n^{1/2}t)dsdt=\w_{2n}(x)+\rho_n(x).
\end{gather*}
Подставляя эти разложения в \eqref{omega}, получим
\begin{multline*}
\omega(x,\la_n)=\sin(nx)\left(1-a_{2n}(x)-\frac1{2n}B_{2n}(x)\right)+\\
+x\cos(nx)\left(-\frac{b_{2n}}{\pi}+\frac{A_{2n}}{2\pi n}-\frac{2\w_{2n}}{\pi}-\frac{U(\pi)}{2\pi n}\right)+\\
+\cos(nx)\left(b_{2n}(x)+\frac{U(x)}{2n}+2\w_{2n}(x)-\frac{A_{2n}(x)}{2n}\right)+\rho_n(x).
\end{multline*}
Остаётся нормировать эти собственные функции. Для этого вычислим $\|\omega(x,\la_n)\|$:
\begin{gather*}
\intl_0^\pi \omega(x,\la_n)\overline{\omega(x,\la_n)}dx=
\intl_0^\pi\sin^2nxdx-2\Re\intl_0^\pi\sin^2(nx)(a_{2n}(x)+\frac1{2n}B_{2n}(x))dx-\\
-\Re\intl_0^\pi\sin(2nx)(-\frac{x}\pi\upsilon_{2n}(\pi)+\upsilon_{2n}(x))dx+\rho_n=\\
=\frac{\pi}2-\intl_0^\pi(\pi-x)\ur(x)\cos(2nx)dx-\frac1{2n}\intl_0^\pi(\pi-x)(\ur^2(x)=\ui^2(x))\sin(2nx)dx+\rho_n.
\end{gather*}
Отсюда сразу же следует формула \eqref{asef}.Теперь найдём асимптотические формулы для векторов биортогональной системы. Поскольку асимптотически все собственные значения просты, имеем $\v_n(x)=\frac{\overline{y_n}}{(y_n,\overline{y_n})}$. Заметим, что
$$
(y_n,\overline{y_n})=\intl_0^\pi y_n^2(x)=1-\frac{2i}\pi\intl_0^\pi(\pi-x)\ui(x)\cos(2nx)dx-
\frac{2i}{\pi n}\intl_0^\pi(\pi-x)\ur(x)\ui(x)\sin(2nx)dx+\rho_n,
$$
откуда получаем формулу \eqref{asbiort}.
\end{proof}

\end{document}